\documentclass[11pt,fleqn]{article}
\usepackage{amsfonts, epsfig, amssymb, amsmath, relsize}
\usepackage[normalem]{ulem}
\usepackage{color}
\newcommand{\ol}{\setlength{\itemsep}{0pt.}\begin{enumerate}}
\newcommand{\eol}{\end{enumerate}\setlength{\itemsep}{-\parsep}}
\newcommand{\ignore}[1]{}
\title{Weight distribution of random linear codes and Krawchouk polynomials}
\author{Alex Samorodnitsky}

\begin{document}
\date{}
\maketitle


\newtheorem{THEOREM}{Theorem}[section]
\newenvironment{theorem}{\begin{THEOREM} \hspace{-.85em} {\bf :}
}%
                        {\end{THEOREM}}
\newtheorem{LEMMA}[THEOREM]{Lemma}
\newenvironment{lemma}{\begin{LEMMA} \hspace{-.85em} {\bf :} }%
                      {\end{LEMMA}}
\newtheorem{COROLLARY}[THEOREM]{Corollary}
\newenvironment{corollary}{\begin{COROLLARY} \hspace{-.85em} {\bf
:} }%
                          {\end{COROLLARY}}
\newtheorem{PROPOSITION}[THEOREM]{Proposition}
\newenvironment{proposition}{\begin{PROPOSITION} \hspace{-.85em}
{\bf :} }%
                            {\end{PROPOSITION}}
\newtheorem{DEFINITION}[THEOREM]{Definition}
\newenvironment{definition}{\begin{DEFINITION} \hspace{-.85em} {\bf
:} \rm}%
                            {\end{DEFINITION}}
\newtheorem{EXAMPLE}[THEOREM]{Example}
\newenvironment{example}{\begin{EXAMPLE} \hspace{-.85em} {\bf :}
\rm}%
                            {\end{EXAMPLE}}
\newtheorem{CONJECTURE}[THEOREM]{Conjecture}
\newenvironment{conjecture}{\begin{CONJECTURE} \hspace{-.85em}
{\bf :} \rm}%
                            {\end{CONJECTURE}}
\newtheorem{MAINCONJECTURE}[THEOREM]{Main Conjecture}
\newenvironment{mainconjecture}{\begin{MAINCONJECTURE} \hspace{-.85em}
{\bf :} \rm}%
                            {\end{MAINCONJECTURE}}
\newtheorem{PROBLEM}[THEOREM]{Problem}
\newenvironment{problem}{\begin{PROBLEM} \hspace{-.85em} {\bf :}
\rm}%
                            {\end{PROBLEM}}
\newtheorem{QUESTION}[THEOREM]{Question}
\newenvironment{question}{\begin{QUESTION} \hspace{-.85em} {\bf :}
\rm}%
                            {\end{QUESTION}}
\newtheorem{REMARK}[THEOREM]{Remark}
\newenvironment{remark}{\begin{REMARK} \hspace{-.85em} {\bf :}
\rm}%
                            {\end{REMARK}}

\newcommand{\thm}{\begin{theorem}}
\newcommand{\lem}{\begin{lemma}}
\newcommand{\pro}{\begin{proposition}}
\newcommand{\dfn}{\begin{definition}}
\newcommand{\rem}{\begin{remark}}
\newcommand{\xam}{\begin{example}}
\newcommand{\cnj}{\begin{conjecture}}
\newcommand{\mcnj}{\begin{mainconjecture}}
\newcommand{\prb}{\begin{problem}}
\newcommand{\que}{\begin{question}}
\newcommand{\cor}{\begin{corollary}}
\newcommand{\prf}{\noindent{\bf Proof:} }
\newcommand{\ethm}{\end{theorem}}
\newcommand{\elem}{\end{lemma}}
\newcommand{\epro}{\end{proposition}}
\newcommand{\edfn}{\bbox\end{definition}}
\newcommand{\erem}{\bbox\end{remark}}
\newcommand{\exam}{\bbox\end{example}}
\newcommand{\ecnj}{\bbox\end{conjecture}}
\newcommand{\emcnj}{\bbox\end{mainconjecture}}
\newcommand{\eprb}{\bbox\end{problem}}
\newcommand{\eque}{\bbox\end{question}}
\newcommand{\ecor}{\end{corollary}}
\newcommand{\eprf}{\bbox}
\newcommand{\beqn}{\begin{equation}}
\newcommand{\eeqn}{\end{equation}}
\newcommand{\wbox}{\mbox{$\sqcap$\llap{$\sqcup$}}}
\newcommand{\bbox}{\vrule height7pt width4pt depth1pt}
\newcommand{\qed}{\bbox}

\newcommand{\bigast}{\mathop{\Large \mathlarger{\mathlarger{*}}}}

\def\sup{^}

\def\H{\{0,1\}^n}

\def\S{S(n,w)}

\def\g{g_{\ast}}
\def\xop{x_{\ast}}
\def\y{y_{\ast}}
\def\z{z_{\ast}}

\def\f{\tilde f}

\def\n{\lfloor \frac n2 \rfloor}

\def \E{\mathop{{}\mathbb E}}
\def \R{\mathbb R}
\def \C{\mathbb C}
\def \Z{\mathbb Z}
\def \F{\mathbb F}
\def \T{\mathbb T}

\def \x{\textcolor{red}{x}}
\def \r{\textcolor{red}{r}}
\def \Rc{\textcolor{red}{R}}

\def \noi{{\noindent}}

\def \iff{~~~~\Leftrightarrow~~~~}
\def \lrarrow{\leftrightarrow}

\def \queq {\quad = \quad}

\def\<{\left<}
\def\>{\right>}
\def \({\left(}
\def \){\right)}

\def \e{\epsilon}
\def \l{\lambda}

\def\Tp{Tchebyshef polynomial}
\def\Tps{TchebysDeto be the maximafine $A(n,d)$ l size of a code with distance $d$hef polynomials}
\newcommand{\rarrow}{\rightarrow}

\newcommand{\larrow}{\leftarrow}

\overfullrule=0pt
\def\setof#1{\lbrace #1 \rbrace}

\begin{abstract}
For $0 < \l < 1$ and $n \rarrow \infty$ pick uniformly at random
$\l n$ vectors in $\H$ and let $C$ be the orthogonal complement
of their span. Given $0 < \gamma < \frac12$ with $0 < \l  < h(\gamma)$,
let $X$ be the random variable that counts the number of words
in $C$ of Hamming weight $i = \gamma n$ (where $i$ is assumed to be an even integer). Linial and Mosheiff \cite{LM} determined the
asymptotics of the moments of $X$ of all orders $o\(\frac{n}{\log n}\)$. In this paper we extend their estimates up to moments of linear order. Our key observation is that the behavior of the suitably normalized $k^{th}$ moment of $X$ is essentially determined by the $k^{th}$ norm of the Krawchouk polynomial $K_i$.

\end{abstract}

\section{Introduction}

\noi This paper follows up on the recent result \cite{LM} by Linial and Mosheiff. In \cite{LM} the authors study the distribution of the number of codewords of a given weight in a random linear code of a given rate, providing tight estimates on the suitably normalized moments of this distribution, up to moments of order $o\(\frac{n}{\log n}\)$. In this paper we extend the estimates in \cite{LM} up to moments of linear order.

\noi We refer to \cite{LM} for the introduction and initial discussion of the problem and for the description of the wider context of asymptotic problems in coding theory. Here we pass directly to technical definitions. In the discussion below we try to adhere to the notation of \cite{LM}, where possible.

\noi {\bf Definitions and notation}.Let $h(t) = t \log_2\(\frac 1t\) + (1-t) \log_2\(\frac{1}{1-t}\)$ be the binary entropy function. Given $0 < \gamma < \frac12$ and $0 < \l < h(\gamma)$, assume that an integer $n$ is such that $\gamma n$ and $\lambda n$ are integers. Furthermore, assume that $i = \gamma n$ is even. Let $L$ be the set of all vectors of weight $i$ in $\H$. Let $M$ be a random binary matrix with $\lambda n$ rows and $n$ columns. Let $C$ be a (random) linear code with parity check matrix $M$. Set $X = |L \cap C|$. Set $p = 2^{-\l n}$.

\noi As observed in \cite{LM}, it is well-known (and easy to see) that $\E X = |L| \cdot \E |C| = {n \choose i} \cdot p$ and $\mathrm{Var}(X) = {n \choose i} \cdot p(1-p) \approx {n \choose i} \cdot p$. Note that here and below we use $\approx$, $\gtrsim$, and $\lesssim$ notation to describe equalities and inequalities which hold up to lower order terms. We also write $\gg$ to indicate that the left hand side of an inequality is exponentially larger than its right hand side.

\noi The following theorem from \cite{LM} describes higher central moments of $X$:

\thm
\label{thm:LM}
Assume $2 \le k \le o\(\frac{n}{\log n}\)$. Then
\[
\frac{(X - \E X)^k}{\mathrm{Var}(X)^{\frac k2}} ~\approx~ \left\{\begin{array}{cc} o(1) & \mbox{if $k$ is odd and $< k_0$} \\ k!! & \mbox{if $k$ is even and $< k_0$} \\ 2^{\(F(k,\gamma) - \frac{k}{2} h(\gamma) - \(\frac{k}{2}-1\) \lambda\) \cdot n} & \mbox{if $k \ge k_0$} \end{array} \right.
\]
\ethm

\noi Here $F(k,\gamma)$ is a certain bivariate function defined in \cite{LM} and $k_0$ is the minimal value of $k$ for which $F(k,\gamma) - \frac{k}{2} h(\gamma) - \(\frac{k}{2}-1\) \lambda > 0$ (it can be seen that this inequality holds for all $k \ge k_0$).

\noi \cite{LM} discusses the question of extending the estimates above to higher values of $k$. In this paper we suggest a different approach to the problem, proving the following result.  

\thm
\label{thm:main}
Let $\e = \min\{\l, h(\gamma) - \l\}$. There exist constants $c_{\e}$ and $C_{\e}$ such that for any $C_{\e} \le k \le c_{\e} n$ holds
\[
\frac{(X - \E X)^k}{\mathrm{Var}(X)^{\frac k2}} ~\approx~ 2^{\(F(k,\gamma) - \frac{k}{2} h(\gamma) - \(\frac{k}{2}-1\) \lambda\) \cdot n}.
\]
\ethm

\noi {\bf Discussion}. As observed in \cite{LM} higher central moments of $X = |C \cap L|$ are larger than those for a general random code $C$ in which each word in chosen from $\H$ uniformly and independently with probability $p$. In that case $X$ is a binomial random variable, $X \sim B\({n \choose i}, p\)$.

\noi Roughly speaking, the $k^{th}$ moment of $X$ counts the number of $k$-tuples in $L$ contained in $C$. It is easy to see that the advantage given by linearity comes from linearly dependent $k$-tuples. Our main observation is that, at least intuitively, for $k$ much smaller than $n$ a "typical" dependent $k$-tuple is a $k$-circuit, that is a linear dependence of the form $x_1 + ... + x_k = 0$, which does not contain any smaller dependencies. In particular, the rank of the set $\left\{x_1,x_2,...,x_k\right\}$  is $k-1$. The number of such dependencies is very close (this will be made more precise below) to $\E \(\sum_{x \in L} w_x\)^k = \E K^k_i$, where $w_x$ is the Walsh-Fourier character  corresponding to $x \in \H$, and $K_i$ is the $i^{th}$ Krawchouk polynomial (all notions will be defined in Section~\ref{sec:prelim} below). The expectation is taken over the uniform distribution on $\H$,  Given the assumption that $i$ is an even integer, we have $\E K_i^k \approx \E |K_i|^k = \|K_i\|_k^k$ (see Section~\ref{sec:prelim} for this as well). As observed in \cite{LM} the probability that a subset $S$ of $\H$ is contained in $C$ is $p^{r(S)}$, where $r(S)$ is the rank of $S$. Hence, if this intuition is correct, the expected number of dependent $k$-tuples in $C$ would be close to $p^{k-1} \|K_i\|_k^k$, and we would expect the advantage given by linearity to be visible iff
\[
p^{k-1} \|K_i\|_k^k ~\gg~ \mathrm{Var}(X)^{\frac k2} ~=~ \({n \choose i} \cdot p\)^{\frac k2} \quad \Longleftrightarrow \quad p^{\frac k2 -1} \|K_i\|_k^k ~\gg~ \({n \choose i}\)^{\frac k2},
\]
in which case we would expect to have (assuming that the advantage provided by linearity is greater than the corresponding binomial moment arising from independent choices)
\[
\frac{(X - \E X)^k}{\mathrm{Var}(X)^{\frac k2}} ~\approx~ \frac{p^{\frac k2 -1} \|K_i\|_k^k}{\({n \choose i}\)^{\frac k2}}.
\]

\noi We would like to understand the exponential behavior of these estimates. Recall that $p = 2^{-\lambda n}$. We have  $|L| = {n \choose i} \approx 2^{h(\gamma) n}$ (by the known estimates on binomial coefficients, see e.g. Theorem 1.4.5. in \cite{van Lint}). We also have (see Section~\ref{sec:prelim}) that $\frac 1n \log_2 \|K_i\|^k_k \approx \psi(k,\gamma) + \frac k2 \cdot h(x)$, where $\psi(k,\gamma)$ is a certain bivariate function defined in \cite{KS2}. Using of all this, and passing to exponents, the constraint $p^{k-1} \|K_i\|_k^k \gg \mathrm{Var}(X)^{\frac k2}$ becomes
\[
\psi(k,\gamma)  > \(\frac k2 -1\) \lambda,
\]
in which case we expect to have
\[
\frac{(X - \E X)^k}{\mathrm{Var}(X)^{\frac k2}} ~\approx~ 2^{\(\psi(k,\gamma) - \(\frac k2 -1\) \lambda\) \cdot n}.
\]

\noi It turns out that these are precisely the results proved in the third claim of Theorem~\ref{thm:LM}~and~in~Theorem~\ref{thm:main}. In fact, the two pertinent bivariate functions defined in \cite{LM} and in \cite{KS2} are the same, up to a multiple of the entropy function:
\[
F(k,\gamma) ~=~ \psi(k,\gamma) + \frac k2 \cdot h(\gamma).
\]
To see that, cf. the definition of $F(k,\gamma)$ and Definition~10 in \cite{LM} with Lemma~2.5 and Proposition~2.6 in \cite{KS2}.

\noi The remainder of the paper is devoted to proving Theorem~\ref{thm:main}. By the preceding discussion, this theorem is an immediate corollary of the two following claims.

\pro
\label{pro:lower bound}
Assume that $k \le \lambda n - 1$. Then
\[
\E_C \(X - \E X\)^k ~\ge~ \frac12 \cdot p^{k -1} \|K_i\|_k^k.
\]
\epro

\noi and

\pro
\label{pro:upper bound}
Let $\e = \min\{\l, h(\gamma) - \l\}$. There exist constants $c_{\e}$ and $C_{\e}$ such that for any $C_{\e} \le k \le c_{\e} n$ holds
\[
\E_C \(X - \E X\)^k ~\lesssim~ p^{k -1} \cdot 2^{\(\psi\(k, \gamma\) + \frac k2 h(\gamma)\) \cdot n}.
\]
\epro

\noi These two claims are proved in Section~\ref{sec:proofs}. Prior to that, some relevant notions and background are provided in Section~\ref{sec:prelim}.

\section{Preliminaries}
\label{sec:prelim}

\noi We collect some relevant facts on Fourier analysis on the boolean cube (see e.g., \cite{O'Donnel}). 

\noi For $x \in \H$, define the Walsh-Fourier character $w_x$ on $\H$ by setting $w_x(y) = (-1)^{\sum_i x_i y_i}$, for all $y \in \H$. The {\it weight} of the character $w_x$ is the Hamming weight $|x|$ of $x$ (that is the number of $1$-coordinates in $x$).  The characters $\{w_x\}_{x \in \H}$ form an orthonormal basis in the space of real-valued functions on $\H$, under the inner product $\<f, g\> = \frac{1}{2^n} \sum_{x \in \H} f(x) g(x)$. The expansion $f = \sum_{x \in \H} \widehat{f}(x) w_x$ defines the Fourier transform $\widehat{f}$ of $f$. We also have the Parseval identity, $\|f\|^2_2 = \sum_{x \in \H} {\widehat{f}}^2(x)$. Here is one additional simple fact that we will require. Let $C$ be a linear subspace of $\H$. Then $\widehat{1_C} = \frac{|C|}{2^n} \cdot 1_{C^{\perp}}$.

\noi {\it Krawchouk polynomials}. For $0 \le i \le n$, let $F_i$ be the sum of all Walsh-Fourier characters of weight $i$, that is $F_i = \sum_{|x| = i} w_x$. It is easy to see that $F_i(x)$ depends only on the Hamming weight $|x|$ of $x$, and it can be viewed as a univariate function on the integer points $0,...,n$, given by the restriction to $\{0,...,n\}$ of the univariate polynomial $K_i = \sum_{k=0}^i (-1)^k {x \choose k} {{n-x} \choose {i-k}}$ of degree $i$. That is, $F_i(x) = K_i(|x|)$. The polynomial $K_i$ is the $i^{th}$ {\it Krawchouk polynomial} (see e.g., \cite{mrrw}). Abusing notation, we will also call $F_i$ the $i^{th}$ Krawchouk polynomial, and write $K_i$ for $F_i$ when the context is clear.
Here are two simple properties of the Krawchouk polynomials which we will need: $K_i(0) = \|K_i\|_2^2 = {n \choose i}$.

\noi {\it Norms}. It is easy to see from the definition of the Walsh-Fourier characters that for any integer $k$ and for any subset $S$ of $\H$ holds $\E \(\sum_{x \in S} w_x\)^k = \Big | \{\(x_1,...,x_k\) \in S^k,~x_1+...+x_k = 0\} \Big |$ (the expectation is taken w.r.t. the uniform probability measure on $\H$) . In particular, if $L$ is the set of all vectors of weight $i$ in $\H$ and $k$ is an even integer, we have
\[
\|K_i\|_k^k ~=~ \E \(\sum_{x \in L} w_x\)^k ~=~ \Big | \{\(x_1,...,x_k\) \in L^k,~x_1+...+x_k = 0\} \Big |.
\]
\noi Let now $k > 2$ be an odd integer. As observed e.g., in Section 2.2 in \cite{KS2}, the $k^{th}$ norm of $K_i$ is essentially attained outside the root region of $K_i$. For an even $i$, $K_i$ is positive outside its root region and hence we also have
\[
\|K_i\|_k^k ~=~ \E |K_i|^k ~\approx~ \E K_i^k ~=~ \Big | \{\(x_1,...,x_k\) \in L^k,~x_1+...+x_k = 0\} \Big |.
\]

\noi As observed e.g., in Section 2.2 in \cite{KS2}, for all $0 \le i \le n/2$ and $k \ge 2$ holds, writing $\gamma = i/n$:

\begin{itemize}

\item $\frac{\|K_i\|_k^k}{\(\|K_i\|_2\)^k} ~\le~ 2^{\psi\(k, \frac in\) \cdot n} = 2^{\psi\(k, \gamma\) \cdot n}$, where $\psi$ is a certain bivariate function defined in \cite{KS2}.

\item Moreover, $\frac{\|K_i\|_k^k}{\(\|K_i\|_2\)^k} ~\approx~ 2^{\psi\(k, \gamma\) \cdot n}$.

\end{itemize}

\noi As pointed out in the introduction, $F(k,\gamma) = \psi(k,\gamma) + \frac k2 \cdot h(\gamma)$. Since $\|K_i\|_2^2 = {n \choose i} \approx 2^{h(\gamma) \cdot n}$, we have
\[
\|K_i\|_k^k ~\approx~ 2^{\(\psi\(k, \gamma\) + \frac k2 h(\gamma)\) \cdot n} ~=~ 2^{F\(k, \gamma\) \cdot n}.
\]

\noi Finally, we need some properties of the function $\psi$ shown in \cite{LM}. We collect them in the following lemma, for convenience (writing them in terms of $\psi$). We remark that the second and the third properties listed in the lemma were also shown independently in \cite{KS2}.

\lem
\label{lem:psi}

\begin{enumerate}

\item $\psi(k, \gamma) \ge \frac k2 \cdot  H(\gamma) - 1$.

\item For a fixed $\gamma$ the function $\psi(k,\gamma)$ is convex in $k$.

\item Since, furthermore, $\psi(2,\gamma) = 0$, this implies that $\frac{\psi(k, \gamma)}{k-2}$ is an increasing function of $k$ for $k > 2$.

\end{enumerate}

\elem

\rem 
It seems worthwhile to sketch a way to derive these facts from the approximate identity 
$\|K_i\|_k^k ~\approx~ 2^{\(\psi\(k, \gamma\) + \frac k2 h(\gamma)\) \cdot n}$. Observe first that $\|K_i\|_k^k \ge \frac{1}{2^n} \cdot K(0)^k = \frac{1}{2^n} \cdot {n \choose i}^k$; and second that for any function $f$ on $\H$, the function $\alpha \rarrow \log \|f\|_{\frac{1}{\alpha}}$ is convex on $(0,1]$ (this is a consequence of H\"older's inequality, see e.g., \cite{HLP}, Theorems 196 and 197).

\erem

\section{Proof of Theorem~\ref{thm:main}}
\label{sec:proofs}

\subsection{Proof of Proposition~\ref{pro:lower bound}}

\noi We will show this proposition to be an immediate corollary of the lower bound given in the following claim. (The upper bound in this claim we be used in the next section as a step towards the proof of Proposition~\ref{pro:upper bound}).

\pro
\label{pro:moment bounds}

Assume that $k \le \lambda n - 1$. Then
\[
\frac 12 \cdot \sum_{S = \(u_1...u_k\) \in L^k} p^{r(S)} ~\le~ \E_C \(X - \E X\)^k ~\le~ 2 \cdot \sum_{S = \(u_1...u_k\) \in L^k} p^{r(S)},
\]
where the summation is over all sequences $\(u_1...u_k\)$ of vectors which contain no coloops (vector which is not contained in the span of all the rest).
\epro

\noi Given the lower bound in Proposition~\ref{pro:moment bounds}, Proposition~\ref{pro:lower bound} can be derived as follows. Note that a sequence of vectors which sums to $0$ necessarily contains no coloops. Hence we have (see Section~\ref{sec:prelim})
\[
p^{k-1} \cdot \|K_i\|_k^k  ~=~ p^{k-1} \cdot \sum_{S = \(u_1...u_k\) \in L^k, u_1+...+u_k = 0} 1 ~\le~
\]
\[
\sum_{S = \(u_1...u_k\) \in L^k, u_1+...+u_k = 0} p^{r(S)} ~\le~  \sum_{S = \(u_1...u_k\) \in L^k} p^{r(S)},
\]
where the last summation is over all sequences $\(u_1...u_k\)$ of vectors which contain no coloops.

\prf (Of Proposition~\ref{pro:moment bounds})

\noi \noi We will use the following notation. For $u \in \H$, let $Y_u$ be the indicator of the event $u$ is in $C$. Let $Z_u = Y_u - p$. Then $X = \sum_{u \in L} Y_u$ and $X - \E X = \sum_{u \in L} Z_u$. We have
\[
\E_C \(X - \E X\)^k ~=~ \E_C \(\sum_{u \in L} Z_u\)^k ~=~ \sum_{\(u_1...u_k\) \in L^k} \E_C \prod_{r=1}^k Z_{u_r}.
\]

\noi So, the claim of the proposition will follow from the next lemma.

\lem
\label{lem:expectation}
Assume that $k \le \lambda n - 1$. Let $S = \(u_1...u_k\)$ be a sequence of vectors in $\H$. Then there are two cases.

\begin{itemize}

\item These vectors contain no coloops. In this case
\[
\frac12 p^{r(S)} ~\le~ \E_C \prod_{r=1}^k Z_{u_r} ~\le~ 2 p^{r(S)}.
\]

\item These vectors contain a coloop. Then
\[
\E_C \prod_{r=1}^k Z_{u_r} ~=~ 0.
\]

\end{itemize}

\elem

\prf (Of Lemma~\ref{lem:matroid})

\noi Let $g = g_S = \prod_{r=1}^k Z_{u_r}$. This is a function on subspaces which depends only on the number of vectors from $S$ contained in the given subspace. We need to estimate $\E g$. It will be convenient for us to transform the setting in the following manner. Let $t$ be the number of distinct vectors in $S$. Arranging them in some order, and writing the statistics of their appearance in $S$, converts $S$ into a $t$-tuple of integers $\(s_1...s_t\)$ summing to $k$. A subspace corresponds to a vector $x \in \{0,1\}^t$ whose $1$-coordinates indicate which distinct vectors from $S$ this subspace contains. The function $g$ then transforms into a function on $\{0,1\}^t$ defined by $g(x) = (1-p)^{\<S,x\>} (-p)^{k - \<S,x\>}$. We have a measure ${\cal L}$ on $\{0,1\}^t$ induced by the measure on subspaces, which is determined by the following property: for all $R \subseteq [t]$ holds ${\cal L}\{x: x_j = 1 ~\forall j \in R\} = p^{r(R)}$, where $r(R)$ stands for $r\(\{u_i\}_{i \in R}\)$.
With all this, we have
\[
\E_C \prod_{r=1}^k Z_{u_r} ~=~ \E_{x \sim {\cal L}} g(x).
\]

\noi Let $f$ be a function on $\{0,1\}^t$ defined for $R \subseteq [t]$ by $f(R) = p^{r(R)}$. Thinking of ${\cal L}$ as a function on $\{0,1\}^t$, we have $f(R) = \sum_{x:\,R \subseteq x} {\cal L}(x)$. Hence, by the M\"obius inversion formula for the boolean lattice (see e.g., \cite{LW}) we have ${\cal L} = \mu f$, where $\mu$ is the M\"obius function of the boolean lattice: $\mu(x,y) = (-1)^{|y|-|x|}$ if $x$ is a subset of $y$, and $0$ otherwise.

\noi Writing $\<\cdot, \cdot\>$ for the inner product of functions on $\{0,1\}^t$ endowed with the counting measure, we have
\[
\E_{z \sim {\cal L}} g(z) ~=~ \<{\cal L},g\> ~=~ \<\mu f, g\> = \<f, \mu^t g\>.
\]

\noi Next, we compute $\mu^t g$. We have
\[
\(\mu^t g\)(x) ~=~ \sum_z \mu^t(x,z) g(z) ~=~ \sum_{z \subseteq x} (-1)^{|x|-|z|} (1-p)^{\<S,z\>} (-p)^{k - \<S,z\>} ~=~
\]
\[
(-1)^{k+|x|} p^k \cdot \sum_{z \subseteq x} (-1)^{\<S,z\>+|z|} \(\frac{1-p}{p}\)^{\<S,z\>}.
\]

\noi Using the fact that $\sum_{z \subseteq x} \prod_{i \in z}\alpha_i = \prod_{i \in x} \(1 + \alpha_i\)$, the last expression is
\[
(-1)^{k+|x|} p^k \prod_{i \in x} \(1 - \(\frac{p-1}{p}\)\)^{s_i}.
\]

\noi Hence
\[
 \<f, \mu^t g\> ~=~ (-1)^{k} p^k  \sum_{x \in \{0,1\}^t} (-1)^{|x|} p^{r(x)} \prod_{i \in x} \(1 - \(\frac{p-1}{p}\)\)^{s_i} ~=~
\]
\[
(-1)^{k} p^k  \sum_{x \in \{0,1\}^t} (-1)^{|x|} p^{r(x) - \<S,x\>} \prod_{i \in x} \(p^{s_i} - (p-1)^{s_i}\).
\]

\noi Now there are two cases to consider.  Let $\sigma_x$ denote the summand corresponding to $x$ in the last expression. Assume first that $S$ contains a coloop (w.l.o.g. vector number $t$). We claim that in this case for all $x$ holds $\sigma_x = -\sigma_{x \oplus e_t}$ and hence the total sum is $0$. In fact, this follows by inspection, since (clearly) $s_t = 1$ and $r(x \cup t) = r(x \setminus t) + 1$.

\noi Assume now that $S$ does not contain coloops. Since by assumption $p = 2^{-\l n} \le 2^{-k-1}$ and since $s_i \ge 1$ for all $i \in [t]$, we have that, up to a $1 + o(1)$ multiplicative factor, $\prod_{i \in x} \(p^{s_i} - (p-1)^{s_i}\) \approx \prod_{i \in x} (-1)^{s_i + 1} = (-1)^{\<S,x\> + |x|}$. So, $\sigma_x \approx (-1)^{k} p^k (-1)^{\<S,x\>} p^{r(x) - \<S,x\>}$. In particular, $\sigma_{[t]} \approx p^{r(S)}$. We claim that this is the dominant term in $\sum_{x \in \{0,1\}^t} \sigma_x$. Specifically, we claim that $\sigma_{[t]} \ge 2 \cdot \sum_{x \not = [t]} |\sigma_x|$. Note that showing this will complete the proof of the lemma.

\noi Since the number of summands is smaller than $\frac{1}{2p}$, it suffices to show that for all $x \not = [t]$ holds $r(x) - \<S,x\> > r(S) - k$. In fact,
\[
(r(x) - \<S,x\>) - (r(S) - k) ~=~ \sum_{i \not \in x} s_i - (r(S) - r(x)).
\]
Since $r(S) - r(x) \le t - |x|$, this can only be $0$ iff $s_i = 1$ for all $i \not \in x$ and $r(S) - r(x) = |t| - |x|$. But this means that all vectors in $x^c$ are coloops, contradicting the assumptions.

\eprf~( Lemma~\ref{lem:matroid}).

\noi This concludes the proof of Proposition~\ref{pro:moment bounds}~and~of~Proposition~\ref{pro:lower bound}.

\subsection{Proof of Proposition~\ref{pro:upper bound}}

\noi Assume $k \le \l n - 1$. From the upper bound of Proposition~\ref{pro:moment bounds}, we have
\[
\E_C \(X - \E X\)^k ~\le~ 2 \cdot \sum_{S = \(u_1...u_k\) \in L^k} p^{r(S)},
\]
where the summation is over all sequences of vectors which contain no coloops. Note that the rank of any such sequence is smaller than $k$. For $1 \le r \le k-1$ let $N(r)$ denote the number of such sequences $S$ with $r(S) = r$. Then the RHS in the last inequality is $2 \cdot \sum_{r=1}^k p^r N(r)$.Observe that $N(k-1)$ is the number of $k$-circuits in $L$. Hence (see Section~\ref{sec:prelim}) $p^{k-1} N(k-1) \le p^{k-1} \cdot \E \|K_i\|_k^k \le p^{k -1} \cdot 2^{\(\psi\(k, \gamma\) + \frac k2 h(\gamma)\) \cdot n}$. Therefore, to prove Proposition~\ref{pro:upper bound} it will suffice to show the following claim.

\pro
\label{pro:upper bound-tech}
Let $\e = \min\{\l, H(\gamma) - \l\}$. There exist constants $c_{\e}$ and $C_{\e}$ such that for any $C_{\e} \le k \le c_{\e} n$ and for any $1 \le r \le k-2$ holds
\[
N(r) ~\lesssim~ p^{k-r-1} \cdot 2^{\(\psi\(k, \gamma\) + \frac k2 h(\gamma)\) \cdot n}.
\]
\epro

\noi In the remainder of this section we prove Proposition~\ref{pro:upper bound-tech}. Let $S = \(u_1...u_k\)$ be a sequence of $k$ vectors in $\H \setminus \{0\}$, with $r(S) = r$. Let $B = \{b_1,...,b_r\}$ be a basis of $S$, that is a subset of $r$ linearly independent vectors in $S$ (which span all vectors in $S$). We build on a simple observation, which we state in the following lemma.

\lem
\label{lem:matroid}
If $S$ contains no coloops, then for some $1 \le v \le \min\{r,k-r\}$ there are additional $v$ vectors $x_1,...,x_v$ in $S$ with the following property: For $1 \le d \le v$, let $x_d = \sum_{i \in S_d} b_i$. Then $S_1,...,S_v$ are non-empty subsets of $[r]$, for all $1 \le d \le v$ holds $S_d \not \subseteq \cup_{i=1}^{d-1} S_i$, and $\cup_{i=1}^v S_i = [r]$.
\elem

\prf
Note that each element of $B$ participates in at least one linear dependence among the elements of $S$. Let $x_1$ be a vector in $S \setminus B$ whose representation as a linear combination of elements of $B$ contains $b_1$. If $S_1 \not = [r]$, let $2 \le j \le r$ be the minimal index in $[r] \setminus S_1$. Let $x_2$ be a vector in $S \setminus B$ whose representation as a linear combination of elements of $B$ contains $b_j$, etc.
\eprf

\noi We proceed with a technical lemma. We write $K$ for $K_i = K_{\gamma n}$ for notational convenience from now on.

\lem
\label{lem:upper}
Let $1 \le v \le r$. Let $m = r + v$. Let $S_1,...,S_v$ be non-empty subsets of $[r]$ with the following property: for all $1 \le d \le v$ holds $S_d \not \subseteq \cup_{i=1}^{d-1} S_i$, and $\cup_{i=1}^v S_i = [r]$.

\noi Let ${\cal A} \subseteq L^m$ contain the $m$-tuples $\(x_1,...,x_m\)$ which satisfy $x_{r+d} = \sum_{i \in S_d} x_i$, $d = 1,...,v$. Then there exist integers $a_1,...,a_v \ge 2$ with $\sum_{i=1}^v a_i = m$ so that
\[
|{\cal A}| ~\leq~ \prod_{i=1}^v \E_{y_i} |K\(y_i\)|^{a_i}.
\]

\noi In fact, the sequence $\{a_d\}$ is defined as follows: $a_d = |S_d \setminus \cup_{i=1}^{d-1} S_i| + 1$, $d = 1,...,v$.
\elem

\prf

\noi Viewing an $m$-tuple of vectors $x_1,...,x_m$ in $\H$ as rows of an $m \times n$ binary matrix $\bar{x}$, let the columns of this matrix be denoted by $c_1\(\bar{x}\),...,c_n\(\bar{x}\)$. Let $M$ be the following $r \times m$ binary matrix: The first $r$ columns of $M$ form an $r \times r$ identity matrix, and for $1 \le d \le v$, the column $r+d$ is the characteristic vector of the set $S_d$. Let $C \subseteq \{0,1\}^m$ be the linear code generated by the rows of $M$. The rows of $M$ are linearly independent, and hence $\dim(C) = r$.

\noi Observe that
\[
{\cal A} ~=~ \big\{\bar{x} = \(x_1,...,x_m\) \in L^m,~c_j\(\bar{x}\) \in C,~j=1,...,n \Big\}
\]

\noi Let
\[
{\cal B} ~=~ \big\{\bar{x} = \(x_1,...,x_m\) \in \(\H\)^m,c_j\(\bar{x}\) \in C,~j=1,...,n \Big\}.
\]
Note that ${\cal A} = {\cal B} \cap L^m$.

\noi Let
\[
{\cal B}^{\ast} ~=~ \big\{\bar{y} = \(y_1,...,y_m\) \in \(\H\)^m,~M c_j\(\bar{y}\) = 0,~j=1,...,n \Big\}
\]

\noi We need an auxiliary lemma.
\lem
\label{lem:upper-aux}
Let $y_1,...,y_m \in \H$. Then
\[
\sum_{\(x_1,...,x_m\) \in {\cal B}} ~\prod_{i=1}^m w_{x_i}\(y_i\) ~=~ \left\{\begin{array}{ccc} 2^{rn} & \mbox{if} & \(y_1,...,y_m\) \in {\cal B}^{\ast} \\ 0 & \mbox{otherwise} \end{array} \right.
\]
\elem

\prf (of Lemma~\ref{lem:upper-aux})

\noi Identifying $\(\H\)^m$ with $\{0,1\}^{mn}$, we view $\bar{x} = \(x_1,...,x_m\), \bar{y} = \(y_1,...,y_m\)$ as points in $\{0,1\}^{mn}$, and then $\prod_{i=1}^m w_{x_i}\(y_i\) = w_{\bar{y}}\(\bar{x}\)$. Note that in this identification ${\cal B}$ becomes an $(rn)$-dimensional linear subspace of $\{0,1\}^{mn}$, and ${\cal B}^{\ast}$ becomes its dual, and we have (see Section~\ref{sec:prelim})
\[
\sum_{\(x_1,...,x_m\) \in {\cal B}} ~\prod_{i=1}^m w_{x_i}\(y_i\) ~=~ \sum_{\bar{x} \in {\cal B}}  w_{\bar{y}}\(\bar{x}\) ~=~ 2^{mn} \cdot \widehat{1_{\cal B}}\(\bar{y}\) ~=
\]
\[
\left\{\begin{array}{ccc} 2^{rn} & \mbox{if} & \bar{y} \in {\cal B}^{\perp} \\ 0 & \mbox{otherwise} \end{array} \right. ~=~
\left\{\begin{array}{ccc} 2^{rn} & \mbox{if} & \(y_1,...,y_m\) \in {\cal B}^{\ast} \\ 0 & \mbox{otherwise} \end{array} \right.
\]

\eprf~(Lemma~\ref{lem:upper-aux})

\noi We continue with the proof of Lemma~\ref{lem:upper}. Since $K = \sum_{z \in L} w_z$ (see Section~\ref{sec:prelim}), for any $x \in \H$ holds $\E_{y \in \H} w_x(y) K(y) = 1_{x \in L}$. Hence, we have, using Lemma~\ref{lem:upper-aux} in the last step,
\[
|{\cal A}| ~=~ \sum_{\(x_1,...,x_m\) \in {\cal B}} ~\prod_{i=1}^m \E_{y_i} w_{x_i}\(y_i\) K_i\(y_i\) ~=~ \sum_{\(x_1,...,x_m\) \in {\cal B}} \E_{y_1,...,y_m} ~\prod_{i=1}^m K\(y_i\) \prod_{i=1}^m w_{x_i}\(y_i\) ~=
\]
\[
\E_{y_1,...,y_m} \prod_{i=1}^m K\(y_i\)  \sum_{\(x_1,...,x_m\) \in {\cal B}} ~\prod_{i=1}^m w_{x_i}\(y_i\) ~=~ 2^{-vn} \sum_{\(y_1,...,y_m\) \in {\cal B}^{\ast}} ~\prod_{i=1}^m K\(y_i\)
\]
Consider the last expression. Let $M^{\ast}$ be the following generating matrix of $C^{\perp}$. It is a $v \times m$ binary matrix whose first $r$ columns form a $v \times r$ matrix whose rows are the characteristic vectors of the sets $S_1,...,S_v$. The remaining $v$ columns form a $v \times v$ identity matrix. Let the subsets of $[v]$ corresponding to the first $r$ columns of $M^{\ast}$ be denoted by $T_1,...,T_r$. Observe that
\[
2^{-vn} \sum_{\(y_1,...,y_m\) \in {\cal B}^{\ast}} ~\prod_{i=1}^m K\(y_i\) ~=~ \E_{y_1,...,y_v \in \H} \prod_{i=1}^v K\(y_i\) \prod_{j = 1}^r K\(\sum_{i \in T_j}y_i\).
\]

\noi For $1 \le d \le v$, set $R_d = S_d \setminus \cup_{i=1}^{d-1} S_i$. Note that the sets $R_1...R_v$ are non-empty and they partition $[r]$.
Note also that the sets $\{T_j\}_{j \in R_d}$ are subsets of $\{d,...,v\}$ and they all contain $d$. We have
\[
\E_{y_1,...,y_v} \prod_{i=1}^v K\(y_i\) \prod_{j = 1}^r K\(\sum_{i \in T_j}y_i\) ~\le~
\]
\[
\E_{y_1,...,y_v} \prod_{i=1}^v |K|\(y_i\) \prod_{j = 1}^r |K|\(\sum_{i \in T_j}y_i\) ~=~
\]
\[
\E_{y_v} |K|\(y_v\) \cdot \prod_{j \in R_v} |K|\(\sum_{i \in T_j}y_i\) \cdot \E_{y_{v-1}} |K|\(y_{v-1}\) \cdot \prod_{j \in R_{v-1}} |K|\(\sum_{i \in T_j}y_i\) \cdots \E_{y_1} |K|\(y_1\) \cdot \prod_{j \in R_1} |K|\(\sum_{i \in T_j}y_i\)
\]

\noi Consider the expectation $E_{y_d} |K|\(y_d\) \cdot \prod_{j \in R_d} |K|\(\sum_{i \in T_j}y_i\)$, for fixed $y_{d+1},...,y_v$. Note that this is the expectation of a product of $a_d = |R_d| + 1 = |S_d \setminus \cup_{i=1}^{d-1} S_i| + 1$ functions $f_1,..., f_{a_d}$ of $y_d$, each of which is a rearrangement of $|K|$. By H\"older's inequality we have
\[
\E_{y_d} f_1 \cdots f_{a_d} ~\le~ \prod_{j=1}^{a_d} \|f_j\|_{a_d} ~=~ \|K\|_{a_d}^{a_d} ~=~ \E_{y_d} |K\(y_d\)|^{a_d}.
\]
Since this holds for all $1 \le d \le v$, we have
\[
\E_{y_v} |K|\(y_v\) \cdot \prod_{j \in R_v} |K|\(\sum_{i \in T_j}y_i\) \cdot \E_{y_{v-1}} |K|\(y_{v-1}\) \cdot \prod_{j \in R_{v-1}} |K|\(\sum_{i \in T_j}y_i\) \cdots \E_{y_1} |K|\(y_1\) \cdot \prod_{j \in R_1} |K|\(\sum_{i \in T_j}y_i\) \le
\]
\[
 \prod_{i=1}^v \E_{y_i} |K\(y_i\)|^{a_i},
\]
completing the proof of the lemma.

\eprf~(Lemma~\ref{lem:upper})

\noi We can now complete the proof of Proposition~\ref{pro:upper bound-tech}. Fix $1 \le r \le k-2$. Let $S = \(u_1,...,u_k\)$ be a sequence of vectors in $L$ of rank $r$ and with no coloops. This sequence is determined fully if we provide the following information: A subset $B$ of $r$ vectors in $S$ which is a basis of $S$, an integer $1 \le v \le \min\{r,k-r\}$, a subset $V$ of $S \setminus B$ of cardinality $v$ which satisfies the conditions of Lemma~\ref{lem:matroid}; and given that, we describe the remaining $k - r - v$ vectors in $S$.

\noi There are ${k \choose r}$ ways to choose the location of the subset $B$, ${{k-r} \choose v} \cdot v!$ ways to describe the location and the ordering of vectors in $V$. Then we need to choose the sets $S_1,...,S_v$ described in Lemma~\ref{lem:matroid}. The number of such sets is at most $2^{rv}$. Provided this information, the number of possible $(r+v)$-tuples of vectors in $B \cup V$ is bounded, by Lemma~\ref{lem:upper}, by $\prod_{i=1}^v \E_{y_i} |K\(y_i\)|^{a_i}$, where $a_1,...,a_v$ are determined by $S_1,...,S_v$. Given such an $(r+v)$-tuple of vectors in $B \cup V$, the remaining $k - r - v$ vectors in $S$ can be chosen in at most $2^{r(k-r-v)}$ ways (since they belong to the span of $B$).

\noi Let $M(r,v) = \max_{\(a_1,...,a_v\)} \prod_{i=1}^v \E_{y_i} |K\(y_i\)|^{a_i}$, where the maximum is taken over all possible integer $v$-tuples $a_1,...,a_v$ which satisfy $a_d \ge 2$ for all $1 \le d \le v$ and $\sum_{d=1}^v a_d = r+v$. Altogether, we get
\[
N(r) ~\le~ \sum_{v=1}^{min\{r,k-r\}} M(r,v) \cdot \min\{r,k-r\} \cdot {k \choose r} \cdot \({{k-r} \choose v} \cdot v!\) \cdot 2^{rv} \cdot 2^{r(k-r-v)} ~\le~
\]
\[
\max_v \big\{M(r,v)\big\} \cdot n \cdot k^{\min\{r,k-r\}} \cdot (k-r)^{\min\{r,k-r\}} \cdot 2^{r(k-r)} ~\le~ \max_v \big\{M(r,v)\big\} \cdot n^{2\min\{r,k-r\} + 1} \cdot 2^{r(k-r)}.
\]
We need to show that, up to lower order terms, this is at most $p^{k-r-1} \cdot 2^{\(\psi\(k, \gamma\) + \frac k2 h(\gamma)\) \cdot n}$. Passing to exponents, it suffices to show that for any $1 \le v \le \min\{r,k-r\}$
\[
\frac 1n \cdot \log_2(M(r,v)) + \frac{2\min\{r,k-r\}+1}{n} \cdot \log_2(n) + \frac{r(k-r)}{n} ~\le~ \psi\(k, \gamma\) + \frac k2 h(\gamma) - (k-r-1) \l.
\]
We refer to Section~\ref{sec:prelim} for facts about the function $\psi$ used in the following argument.

\noi First, we upperbound $M(r,v)$. Let $a_1,...,a_v$ satisfy $a_d \ge 2$ for all $1 \le d \le v$ and $\sum_{d=1}^v a_d = r+v$. Then
\[
\frac 1n \log_2\(\prod_{i=1}^v \E_{y_i} |K\(y_i\)|^{a_i}\) ~\le~ \sum_{i=1}^v \Big(\psi\(a_i,\gamma\) + \frac{a_i}{2} h(\gamma)\Big) ~=~  \sum_{i=1}^v \psi\(a_i,\gamma\) + \frac{r+v}{2} h(\gamma).
\]

\noi Using the fact that $\frac{\psi(a,\gamma)}{a-2}$ increases in $a$ for $a > 2$, we have that $\sum_{i=1}^v \psi\(a_i,\gamma\) \le \frac{r-v}{k-2} \cdot \psi(k,\gamma)$. Hence $\frac 1n \log_2(M(r,v)) \le \frac{r-v}{k-2} \cdot \psi(k,\gamma) + \frac{r+v}{2} h(\gamma)$. Therefore, it suffices to show that
\[
\frac{2\min\{r,k-r\}+1}{n} \cdot \log_2(n) + \frac{r(k-r)}{n} ~\le~ \frac{k-r+v-2}{k-2} \cdot \psi(k,\gamma) + \frac{k-r-v}{2} \cdot h(\gamma) - (k-r-1) \l.
\]

\noi We now use the fact that $\psi(k,\gamma) \ge \frac{k}{2} \cdot h(\gamma) - 1 = \frac{k-2}{2} \cdot h(\gamma) - (1 - h(\gamma))$. Substituting in the above expression, we need to show that
\[
\frac{2\min\{r,k-r\}+1}{n} \cdot \log_2(n) + \frac{r(k-r)}{n} +  \frac{k-r+v-2}{k-2} \cdot (1 - h(\gamma)) ~\le~ (k-r-1) \cdot \(h(\gamma) - \l\).
\]

\noi Recall that $h(\gamma) - \l \ge \e$, for a positive absolute constant $\e$. It is easy to see that there exist constants $c_{\e}$ and $C_{\e}$ such that for sufficiently large $n$ and $C_{\e} \le k \le c_{\e} n$ each of the three summands on the LHS is upperbounded by $\frac{k-r-1}{3}\cdot \e$, completing the proof of the proposition.

\eprf

\end{document}